\documentclass[12pt,leqno,english]{article}
\usepackage{amsmath}
\usepackage{amsfonts}
\usepackage{amssymb}

\newcommand{\R}{\mathbb{R}}
\newcommand{\Z}{\mathbb{Z}}

\newcommand{\rg}{\rightarrow}

\begin{document}

\title
{Winding numbers and Fourier series}

\author{Jean--Pierre Kahane}

\date{}

\maketitle

This is an expository talk on a topic of classical analysis, arizing from the $BMO$--theory of topological degree (Br\'ezis--Nirenberg 1995) \cite{brni}. We sketch the history of the subject and some of its recent developments.

\section{The starting point}

In October 1995 Ha\"{i}m Br\'ezis visited Rutgers University and gave a lecture at the seminar of Israel Gelfand on his recent work with Louis Nirenberg on the extension of the notion of topological degree to a class of functions larger than the class of continuous functions, namely $VMO$, the class of functions with vanishing mean oscillation. Vanishing mean oscillation is expressed by the formula
$$
\lim_{|B|\rg 0} \dfrac{1}{|B|^2} \int\!\!\!\int_{B\times B} |f(x)-f(y)| \ dx\,dy = 0\,
$$
when $f$ maps an open set $\Omega \subset \R^m$ into $\R^n$, $B$ denotes a ball contained in $\Omega$ and $|B|$ its volume, $dx=dx_1\cdots dx_m$ and $|z| = |z_1^2 +\cdots+z_n^2|^{1/2}$. When $f$ maps a manifold into a manifold, this definition should be translated accordingly.

Gelfand was eager to have examples, and Br\'ezis had plenty of them. First,
$$
W^{s,p}(\Omega) \subset VMO(\Omega)
$$
when $0<s<1$, $1<p<\infty$, $sp=n$, and $f\in W^{s,p}(\Omega)$ means
$$
\int\!\!\!\int_{\Omega\times\Omega} \dfrac{|f(x)-f(y)|^p}{|x-y|^{n+sp}} \ dxdy <\infty\,.
$$
Then, the Sobolev classes
$$
H^{n/2} (\Omega) \subset VMO(\Omega)
$$
since $H^s=W^{s,2}$. In particular, going to mapping of the $n$--sphere into itself,
$$
H^1(S^2,S^2) \subset VMO (S^2)\,,
$$
the original motivation of Br\'ezis \cite{brco}. In the same way,
$$
H^{1/2}(S^1,S^1) \subset VMO (S^1)\,.
$$

Gelfand felt satisfied with $S^1$, but not with the definition of $H^{1/2}$ by means of a double integral. Since $f$ can be expressed as $f(e^{it})$, $t\in \R$, and $|f(e^{it})|=1$, how to write the definition in terms of the Fourier coefficients of~$f$,
$$
a_n = \int f(e^{it})e^{-int} \dfrac{dt}{2\pi}
$$
(here and in the sequel we write $\int$ instead of $\int_0^{2\pi}$). Br\'ezis had the answer immediately : the definition can be expressed~as
$$
\sum_{-\infty}^\infty |n|\,|a_n|^2 < \infty\,.
$$

Then Gelfand asked another question. Can you express the topological degree by means of the $a_n$ ? Br\'ezis could not answer on the spot. He went home, made a little computation, and got the simple and beautiful formula
$$
\deg f= \sum_{-\infty}^\infty n|a_n|^2\,.
$$
Here $\deg f$ denotes the usual topological degree if $f$ is continuous, and the $VMO$--degree in general, that is, $f\in H^{1/2}$.

This was integrated in the article of Br\'ezis and Nirenberg \cite{brni} and initiated a series of other questions. The first is already in \cite{brni}. Most of them can be found in the ``mise au point'' by Br\'ezis in 2006 \cite{bre1}. We shall see some answers and new questions.

\vskip2mm

\textbf{Question 1} \cite{brni}.

What happens when $f$ is continuous (then $\deg f$ exists) and does not belong to $H^{1/2}(S^1,S^1)$, that is
$$
\sum_{-\infty}^\infty |n|\,|a_n|^2 =\infty \ ?
$$
Is there any summation process for the series
$$
\sum_{-\infty}^\infty n\,|a_n|^2
$$
such that $\deg f$ can be computed in that way ?

The first answers were given by Jacob Korevaar in 1999 \cite{kor}. Korevaar considers two summation processes, namely
$$
\lim_{n\rg \infty} \sum_{-n}^n m|a_m|^2\,,
$$
using symmetrical partial sums, and
$$
\lim_{r\nearrow 1} \sum_{-\infty}^\infty r^{|m|} |a_m|^2\,,
$$
the process of Abel--Poisson. The second is stronger than the first. Korevaar shows that they work when $f$ has bounded variation, $f\in C\cap BV$, but that none of them work under the mere assumption $f\in C$. Actually they diverge for some $f$, and they converge to a value different from $\deg f$ for some other~$f$.

That led to another question.

\vskip2mm

\textbf{Question 2} \cite{bre1}.

Does $\deg f$ depend on the absolute value of the Fourier coefficients $a_n$ only ? Since the energy of $f$ is defined by the $|a_n|^2$, and the topological degree of a mapping from $S^1$ to $S^1$ is nothing but the winding number, the question can be asked in a pleasant form : can we hear the winding number~?

The answer is negative. Jean Bourgain and Gaby Kozma were able to construct two functions belonging to $C(S^1,S^1)$ with the same $|a_n|$ and different degrees \cite{boko}. It is a very difficult construction.

\vskip2mm

\textbf{Question 3} \cite{kah,kah1}.

Let us return to the summation processes. Let us begin by the Abel--Poisson process, since it is stronger than most usual processes. If we assume that $f$ satisfies a H\"older condition of order $\alpha>0$, that is, in Zygmund's notation, $f\in \Lambda_\alpha(S^1,S^1)$ (Zygmund, as many authors, says ``Lipschitz condition of order $\alpha$'') \cite{zyg}, is it true that
$$
\deg f= \lim_{r\uparrow 1} \sum_{-\infty}^\infty n\ r^{|n|} |a_n|^2\ ?
$$

The answer is positive for $\alpha>\frac{1}{3}$ and negative for $\alpha \le \frac{1}{3}$. This comes from a more precise statement, which involves the classes $\lambda_\alpha^p$ of Zygmund (\cite{zyg}, p.~45), defined as
$$
\lambda_\alpha^p = \{g : \int |g(t+s) -g(s)|^p ds=o(t^\alpha),\ t \downarrow 0\}\,,
$$
and a non--classical summation process, namely
$$
\lim_{t\rg 0}\sum_{-\infty}^\infty |a_n|^2 \dfrac{\sin nt}{t}\,.
$$
This limit exists and is equal to $\deg f$ when $f\in C\cap \lambda_{1/3}^3 (S^1,S^1)$ but there exists $f\in \Lambda_{1/3}(S_1,S_1)$ such that it doesn't exist, or it exists and is different from $\deg f$~\cite{kah}.

The positive part of this statement is valid when $C\cap \lambda_{1/3}^3(S^1,S^1)$ is replaced by $W_{1/3}^3 (S^1,S^1)$ \cite{bre1}. We don't know if it is still valid under the assumption $f\in VMO \cap \lambda_{1/3}^3(S^1,S^1)$ ; that would provide a commun generalization to both statements.

\vskip 2mm

\textbf{Question 4} \cite{kah1}.

We introduced three summation processes, and there are many others. It was a popular subject in the 1920's, and the best reference is Hardy's book of 1949, Divergent Series \cite{har}. I returned to this topic in~\cite{kah1}.

Hardy considers series of terms indexed by positive integers, say
$$
\sum_1^\infty u_m\,.
$$
In our situation
$$
u_m = m (|a_m|^2 -|a_{-m}|^2)\,.
$$
Ordinary convergence to $s$ means
$$
(C)\qquad s= \lim_{n\rg \infty} \sum_1^n u_m
$$
Ces\`aro summability of order $k$ $(k>-1)$ to $s$ means
$$
(C,k)\qquad s=\lim_{n\rg \infty} \Big(
\begin{array}{c}
n + k\\ k
\end{array}\Big)^{-1} 
\sum_m 
\Big(\begin{array}{c}
n + k -m\\
k\end{array}\Big) u_m\, ;
$$
$(C,O)$ is the same as $C$, and $(C,1)$ deals with the arithmetic means of partial sums ; it is the process used  by Fej\'er in Fourier series. The processes
$$
(R,k) \qquad s = \lim_{r\downarrow 0} \sum_1^\infty u_m\Big( \dfrac{\sin mt}{mt}\Big)^k
$$
($k$ being a positive integer) are called Riemann summation processes of order $k$. The original Riemann process is $(R,2)$ and it was used in order to study everywhere convergent trigonometric series. The process we just used is $(R,1)$. Whenever we consider $(R,k)$ we assume $\sum\limits_1^\infty |u_m| m^{-k}< \infty$. The process
$$
(A) \qquad s=\lim_{n\uparrow 1} \sum_1^\infty r^m u_m
$$
is the Abel, or Abel--Poisson, summation process.

It is classical, and easy to see, that $(C,k')$ is stronger than $(C,k)$ and $(R,k')$ is stronger than $(R,k)$ if $k'>k$, and that $(A)$ is stronger than all $(C,k)$ and stronger than $(R,2)$. We write
$$
\begin{array}{lllll}
k'>k &\Longrightarrow &(C,k) &\longrightarrow &(C,k')\\
k'>k &\Longrightarrow &(R,k) &\longrightarrow &(R,k')
\end{array}
$$
$$
\begin{array}{lll}
(C,k) &\longrightarrow &(A)\\
(R,2) &\longrightarrow &(A)
\end{array}
$$
Hardy mentioned some more subtle results :
$$
\begin{array}{llll}
(R,2) &\longrightarrow &(C,2+\delta) &(\delta>0)\\
(R,1) &\longrightarrow &(C,1+\delta) &(\delta>0)
\end{array}
$$
The first is due to Kuttner (1935) \cite{kut} and the second to Zygmund (1928) \cite{zyg2}. Moreover, Kuttner gave an example showing that
$$
(R,3) \nrightarrow (A)\,.
$$

It is easy to see that
$$
(R,1) \nrightarrow (C)\,.
$$
More interesting is the fact that
$$
(R,1) \nrightarrow (C,1)\,.
$$
In the opposite direction,
$$
(C) \nrightarrow (R,1)\,.
$$
These last results can be found at the end of \cite{kah1}.

\vskip2mm

\textbf{Question 5} \cite{bre2,boka}.

In September 2008, Br\'ezis made a report on topological degree and asked the following question.
Let $f\in C(S^1,S^1)$, with Fourier coefficients $a_n(n\in \Z)$. Is it true that
$$
\sum_{-\infty}^\infty |n|\,|a_n|^2 \le |\deg f| + 2 \sum_1^\infty n|a_n|^2 \, ?
$$
Since it is true when the first member is finite, the question can be written~as
$$
\sum_{1}^\infty n\,|a_n|^2 < \infty \Longrightarrow  \sum_{-\infty}^\infty |n|\,|a_n|^2 < \infty\, ?
$$
I could answer in a particular case $(f\in \Lambda_\alpha$, $\alpha>0)$ and Bourgain in the general case.

The answer is positive. Moreover
$$
\sum_{-\infty}^\infty |n|\,|a_n|^2 \le 32 \sum_{1}^\infty n\,|a_n|^2 \,.
$$

The implication is valid in a more general situation. Let $s>0$. Then
$$
\sum_{1}^\infty n^{2s}\,|a_n|^2  < \infty \Longrightarrow \sum_{-\infty}^\infty |n|^{2s}\,|a_n|^2  < \infty\,.
$$
But, except when $s=\dfrac{1}{2}$, there is no constant $C=C(s)$ such that
$$
\sum_{-\infty}^\infty |n|^{2s}\,|a_n|^2  < C \sum_{1}^\infty n^{2s}\,|a_n|^2\qquad \cite{boka}
$$

\vskip2mm

These results are valid when $f$ is supposed to be $VMO$ instead of continuous, but not when $f$ is supposed to be bounded (counterexample : a Blaschke product).

The study of Fourier series of continuous or $VMO$ unimodular functions is an interesting byproduct of the study of winding numbers.

\end{document}